\documentclass{amsart}
\usepackage{amsmath,amssymb,amsfonts}
\usepackage[mathscr]{eucal}

\usepackage{enumerate}
\newenvironment{enumroman}{\begin{enumerate}[\upshape (i)]}
                                                {\end{enumerate}}

\theoremstyle{plain}
\newtheorem{theorem}{Theorem}[section]

\newtheorem{prop}[theorem]{Proposition}

\theoremstyle{definition}
\newtheorem{definition}[theorem]{Definition}
\newtheorem{example}[theorem]{Example}
\newtheorem*{thank}{Acknowledgments}

\input xy
\xyoption{all}

\newcommand{\Deltaop}{{\bf \Delta}^{op}}

\newcommand{\hocolim}{\text{hocolim}}
\newcommand{\holim}{\text{holim}}
\newcommand{\nerve}{\text{nerve}}
\newcommand{\iso}{\text{iso}}
\newcommand{\we}{\text{we}}
\newcommand{\diag}{\text{diag}}
\newcommand{\Hom}{\text{Hom}}
\newcommand{\Map}{\text{Map}}
\newcommand{\Aut}{\text{Aut}}
\newcommand{\Ho}{\text{Ho}}
\newcommand{\SSets}{\mathcal{SS}ets}

\newcommand{\Sets}{\mathcal Sets}
\newcommand{\Secat}{\mathcal Se \mathcal Cat}
\newcommand{\map}{\text{map}}

\newcommand{\hoequiv}{\text{hoequiv}}
\newcommand{\sco}{\mathcal{SC}_\mathcal O}
\newcommand{\css}{\mathcal{CSS}}
\newcommand{\sesp}{\mathcal Se \mathcal Sp}

\begin{document}

\title[Complete Segal spaces]{Complete Segal spaces arising from simplicial categories}

\author[J.E. Bergner]{Julia E. Bergner}

\address{Kansas State University \\ 138 Cardwell Hall
Manhattan, KS 66506}

\email{bergnerj@member.ams.org}

\date{\today}



\begin{abstract}
In this paper, we compare several functors which take simplicial
categories or model categories to complete Segal spaces, which are
particularly nice simplicial spaces which, like simplicial
categories, can be considered to be models for homotopy theories.
We then give a characterization, up to weak equivalence, of
complete Segal spaces arising from these functors.
\end{abstract}

\maketitle

\section{Introduction and Overview}

The idea that simplicial categories, or categories enriched over
simplicial sets, model homotopy theories goes back to a series of
several papers by Dwyer and Kan \cite{dkcalc}, \cite{dkdiag},
\cite{dkfncxes}, \cite{dksimploc}. Taking the viewpoint that a
model category, or more generally a category with weak
equivalences, can be considered to be a model for a homotopy
theory, they develop two methods to obtain from a model category a
simplicial category.  This ``simplicial localization" encodes
higher-order structure which is lost when we pass to the homotopy
category associated to the model category. Furthermore, they prove
that, up to a natural notion of weak equivalence of simplicial
categories (see Definition \ref{dkequiv}), every simplicial
category arises as the simplicial localization of some category
with weak equivalences \cite[2.1]{dkdiag}. Thus, the category of
all (small) simplicial categories with these weak equivalences can
be regarded as the ``homotopy theory of homotopy theories."

This notion was mentioned briefly at the end of Dwyer and
Spalinski's introduction to model categories \cite[11.6]{ds}, and
was made precise by the author in \cite[1.1]{simpcat}, in which
the category $\mathcal{SC}$ of small simplicial categories with
these weak equivalences is shown to have the structure of a model
category.

However, this category is not practically useful for many
purposes.  Simplicial categories are not particularly easy objects
to work with, and the weak equivalences, while natural
generalizations of equivalences of categories, are difficult to
identify.  Motivated by this problem, Rezk defines a model
structure, which we denote $\mathcal{CSS}$, on the category of
simplicial spaces, in which the fibrant-cofibrant objects are
called complete Segal spaces (see Definition \ref{css}). This
model structure is especially nice, in that the objects are just
diagrams of simplicial sets and the weak equivalences, at least
between complete Segal spaces, are just levelwise weak
equivalences of simplicial sets.  Furthermore, this model
structure has the additional structures of a simplicial model
category and a monoidal model category. Because it is given by a
localization of a model structure on the category of simplicial
spaces with levelwise weak equivalences, it is an example of a
presentation for a homotopy theory as described by Dugger
\cite{dugger}.

In his paper, Rezk defines two different functors, one from the
category of simplicial categories, and one from the category of
model categories, to the category of complete Segal spaces. Thus,
he describes a relationship between simplicial categories and
complete Segal spaces, but he does not give an inverse
construction.  However, the author was able to show in
\cite{thesis} that the model categories $\mathcal{SC}$ and
$\mathcal{CSS}$ are Quillen equivalent to one another, and
therefore are models for the same homotopy theory.  Thus, the hope
is that we can answer questions about simplicial categories by
working with complete Segal spaces. This paper is the beginning of
that project.

However, the functor that we use to show that the two model
categories are Quillen equivalent is not the same as Rezk's
functor.  Furthermore, the question arises whether Rezk's functor
on model categories agrees with the composite of the simplicial
localization functor with his functor on simplicial categories.
Rezk gives the beginning of a proof that the two functors agree
when the model category in question has the additional structure
of a simplicial model category.  Our goal in this paper is prove
that, up to weak equivalence in $\mathcal{CSS}$, the two functors
from $\mathcal{SC}$ to $\mathcal{CSS}$ are the same, that Rezk's
result holds for model categories which are not necessarily
simplicial, and that this result does imply that the functor on
model categories agrees with the one on their corresponding
simplicial categories up to weak equivalence.

We then go on to characterize, up to weak equivalence, the
complete Segal space arising from any simplicial category.  It can
be described at each level as the nerve of the monoid of self weak
equivalences of representing objects of the category.

We should mention that the model categories $\mathcal {SC}$ and
$\mathcal{CSS}$ are only two of several known models for the
homotopy theory of homotopy theories.  Our proof that the two are
Quillen equivalent actually uses two intermediate model structures
$\Secat_c$ and $\Secat_f$ on the category of Segal precategories,
or simplicial spaces with a discrete Space at level zero
\cite{thesis}.  The fibrant-cofibrant objects in these structures
are known as Segal categories.  Furthermore, Joyal and Tierney
have shown that there is a Quillen equivalence between each of
these model structures and a model structure $\mathcal{QC}at$ on
the category of simplicial sets \cite{joyalsc}, \cite{jt}. The
fibrant-cofibrant objects of this model structure are known as
quasi-categories and are generalizations of Kan complexes
\cite{joyalquasi}.  One can compare our composite functor
$\mathcal{SC} \rightarrow \css$ to one they define which factors
through $\mathcal{QC}at$ rather than through $\Secat_c$ and
$\Secat_f$.  It is a consequence of Joyal and Tierney's work that
these two different functors give rise to weakly equivalent
complete Segal spaces, as we will describe further in the section
on complete Segal spaces.  An introduction to each of these model
structures and the Quillen equivalences can be found in the survey
paper \cite{survey}.

In \cite{fiberprod}, we use the results of this paper to consider
the complete Segal spaces arising from the homotopy fiber product
construction for model categories, as described by To\"en in his
work on derived Hall algebras \cite{toendha}.  It seems that
results relating diagrams of model categories to diagrams of
complete Segal spaces will prove to be useful.

\begin{thank}
I would like to thank Bill Dwyer, Charles Rezk, Andr\'e Joyal, and
Myles Tierney for conversations about the work in this paper, as
well as the referee for helpful comments on the exposition.
\end{thank}

\section{Background on Model Categories and Simplicial Objects}

In this section, we summarize some facts about model categories,
simplicial sets, and other simplicial objects that we need in the
course of this paper.

A \emph{model category} $\mathcal M$ is a category with three
distinguished classes of morphisms, fibrations, cofibrations, and
weak equivalences.  A morphism which is both a (co)fibration and a
weak equivalence is called an \emph{acyclic (co)fibration}.  The
category $\mathcal M$ with these choices of classes is required to
satisfy five axioms \cite[3.3]{ds}.
%
%
%
%
%
%

Axiom MC1 guarantees that $\mathcal M$ has small limits and
colimits, so in particular $\mathcal M$ has an initial object and
a terminal object.  An object $X$ in a model category $\mathcal M$
is \emph{fibrant} if the unique map $X \rightarrow \ast$ to the
terminal object is a fibration. Dually, $X$ is \emph{cofibrant} if
the unique map from the initial object $\phi \rightarrow X$ is a
cofibration.

The factorization axiom (MC5) can be applied in such a way that,
given any object $X$ of $\mathcal M$, we can factor the map $X
\rightarrow \ast$ as a composite
\[ \xymatrix@1{X \ar[r]^\sim &  X' \ar@{->}[r] & \ast} \] of an
acyclic cofibration followed by a fibration. In this case, $X'$ is
called a \emph{fibrant replacement} of $X$, since it is weakly
equivalent to $X$ and fibrant.  This replacement is not
necessarily unique, but in all the model categories we consider
here it can be assumed to be functorial \cite[1.1.1]{hovey}.
Cofibrant replacements can be defined dually.

The structure of a model category enables us to invert the weak
equivalences formally in such a way that we still have a set,
rather than a proper class, of morphisms between any two objects.
If we were merely to take the localization $\mathcal W^{-1}
\mathcal M$, there would be no guarantee that we would not have a
proper class.  However, with the structure of a model category, we
can define the \emph{homotopy category} $\Ho(\mathcal M)$ to have
the same objects as $\mathcal M$, and as morphisms
\[ \Hom_{\Ho (\mathcal M)}(X,Y) = [X^{cf},Y^{cf}]_\mathcal M, \]
where the right-hand side denotes homotopy classes of maps between
fibrant-cofibrant replacements of $X$ and $Y$, respectively.


The standard notion of equivalence of model categories is given by
the following definitions.  First, recall that an \emph{adjoint
pair} of functors $F \colon \mathcal C \leftrightarrows \mathcal D
\colon G$ satisfies the property that, for any objects $X$ of
$\mathcal C$ and $Y$ of $\mathcal D$, there is a natural
isomorphism
\[ \varphi: \Hom_\mathcal D(FX, Y) \rightarrow \Hom_\mathcal C(X,
GY). \]  The functor $F$ is called the \emph{left adjoint} and $G$
the \emph{right adjoint} \cite[IV.1]{macl}.

\begin{definition} \cite[1.3.1]{hovey}
An adjoint pair of functors $F \colon \mathcal M \leftrightarrows
\mathcal N \colon G$ between model categories is a \emph{Quillen
pair} if $F$ preserves cofibrations and $G$ preserves fibrations.
\end{definition}

\begin{definition} \cite[1.3.12]{hovey}
A Quillen pair of model categories is a \emph{Quillen equivalence}
if if for all cofibrant $X$ in $\mathcal M$ and fibrant $Y$ in
$\mathcal N$, a map $f \colon FX \rightarrow Y$ is a weak
equivalence in $\mathcal D$ if and only if the map $\varphi f
\colon X \rightarrow GY$ is a weak equivalence in $\mathcal M$.
\end{definition}

An important example of a model category is that of the standard
model structure on the category of simplicial sets $\SSets$.
Recall that a \emph{simplicial set} is a functor $X \colon
\Deltaop \rightarrow \Sets$, where $\Deltaop$ is the opposite of
the category ${\bf \Delta}$ of finite ordered sets $[n] = \{0
\rightarrow 1 \rightarrow \cdots \rightarrow n\}$ and
order-preserving maps between them.  We denote the set $X([n])$ by
$X_n$.  In particular in $X$ we have face maps $d_i:X_n
\rightarrow X_{n-1}$ and degeneracy maps $s_i:X_n \rightarrow
X_{n+1}$ for each $0 \leq i \leq n$, satisfying several
compatibility conditions.  Three particularly useful examples of
simplicial sets are the $n$-simplex $\Delta [n]$ and its boundary
$\dot \Delta [n]$ for each $n \geq 0$, and the boundary with the
$k$th face removed, $V[n,k]$, for each $n\geq 1$ and $0 \leq k
\leq n$.  Given a simplicial set $X$, we can take its
\emph{geometric realization} $|X|$, which is a topological space
\cite[I.1]{gj}.

In the standard model category structure on $\SSets$, the weak
equivalences are the maps $f \colon X \rightarrow Y$ for which the
geometric realization $|f|\colon |X| \rightarrow |Y|$ is a weak
homotopy equivalence of topological spaces \cite[I.11.3]{gj}.  In
fact, this model structure on simplicial sets is Quillen
equivalent to the standard model structure on the category of
topological spaces \cite[3.6.7]{hovey}.

More generally, a \emph{simplicial object} in a category $\mathcal
C$ is a functor $\Deltaop \rightarrow \mathcal C$.  The two main
examples which we consider in this paper are those of simplicial
spaces (also called bisimplicial sets), or functors $\Deltaop
\rightarrow \SSets$, and simplicial groups.  We denote the
category of simplicial spaces by $\SSets^{\Deltaop}$.

A simplicial set $X$ can be regarded as a simplicial space in two
ways.  It can be considered a constant simplicial space with the
simplicial set $X$ at each level, and in this case we will denote
the constant simplicial set by $cX$ or just $X$ if no confusion
will arise.  Alternatively, we can take the simplicial space,
which we denote $X^t$, for which $(X^t)_n$ is the discrete
simplicial set $X_n$.  The superscript $t$ is meant to suggest
that this simplicial space is the ``transpose" of the constant
simplicial space.

A natural choice for the weak equivalences in the category
$\SSets^{\Deltaop}$ is the class of levelwise weak equivalences of
simplicial sets.  If we define the cofibrations to be levelwise
also, we obtain a model structure which is usually referred to as
the \emph{Reedy model structure} on $\SSets^{\Deltaop}$
\cite{reedy}.

The Reedy model structure has the additional structure of a
simplicial model category.  In particular, given any two objects
$X$ and $Y$ of $\SSets^{\Deltaop}$, there is a mapping space, or
simplicial set $\Map(X,Y)$ satisfying compatibility conditions
\cite[9.1.6]{hirsch}. In the case where $X$ is cofibrant (as is
true of all objects in the Reedy model structure) and $Y$ is
fibrant, this choice of mapping space is homotopy invariant.

One way to obtain other model structures on the category of
simplicial spaces is to localize the Reedy structure with respect
to a set of maps. While this process works for much more general
model categories \cite[3.3.1]{hirsch}, we will focus here on this
particular case. Let $S=\{f \colon A \rightarrow B\}$ be a set of
maps of simplicial spaces.  A Reedy fibrant simplicial space $W$
is $S$-\emph{local} if for each map $f \in S$, the induced map
\[ \Map(f,W) \colon \Map(B,W) \rightarrow \Map(A,W) \] is a weak
equivalence of simplicial sets.  A map $g \colon X \rightarrow Y$
is then an $S$-\emph{local equivalence} if for any $S$-local
object $W$, the induced map
\[ \Map(g,W) \colon \Map(Y,W) \rightarrow \Map(X,W) \] is a weak
equivalence of simplicial sets.

\begin{theorem} \cite[4.1.1]{hirsch}
There is a model structure $L_S\SSets^{\Deltaop}$ on the category
of simplicial spaces in which
\begin{itemize}
\item the weak equivalences are the $S$-local equivalences,

\item the cofibrations are levelwise cofibrations of simplicial
sets, and

\item the fibrant objects are the $S$-local objects.
\end{itemize}
Furthermore, this model category has the additional structure of a
simplicial model category.
\end{theorem}

We now turn to a few facts about simplicial groups, or functors
from $\Deltaop$ to the category of groups.  Given a simplicial
group $G$, we can take its nerve, a simplicial space with
\[ \nerve(G)_{n,m}= \Hom([m], G_n). \]  Taking the diagonal of
this simplicial space, we obtain a simplicial set, also often
called the nerve of $G$.

From another perspective, for $G$ a simplicial group (or, more
generally, a simplicial monoid), we can find a classifying complex
of $G$, a simplicial set whose geometric realization is the
classifying space $BG$. A precise construction can be made for
this classifying space by the $\overline W G$ construction
\cite[V.4.4]{gj}, \cite{may}. However, we are not so concerned
here with the precise construction as with the fact that such a
classifying space exists, so for simplicity we will simply write
$BG$ for the classifying complex of $G$.

\section{Simplicial Categories and Simplicial Localizations}

In this section, we consider simplicial categories and show how
they arise from Dwyer and Kan's simplicial localization
techniques.  We then discuss model category structures, first on
the category of simplicial categories with a fixed object set, and
then on the category of all small simplicial categories.

First of all, we clarify some terminology. In this paper, by
(small) ``simplicial category" we will mean a category with a set
of objects and a simplicial set of morphisms $\Map(x,y)$ between
any two objects $x$ and $y$, also known as a category enriched
over simplicial sets.  This notion does not coincide with the more
general one of a simplicial object in the category of small
categories, in which we would also have a simplicial set of
objects.  Using this more general definition, if we impose the
additional condition that all face and degeneracy maps are the
identity on the objects, then we get our more restricted notion.

A simplicial category can be seen as a generalization of a
category, since any ordinary category can be regarded as a
simplicial category with a discrete mapping space.  Given any
simplicial category $\mathcal C$, we can consider its
\emph{category of components} $\pi_0 \mathcal C$, which has the
same objects as $\mathcal C$ and whose morphisms are given by
\[ \Hom_{\pi_0 \mathcal C}(x,y) = \pi_0 \Map_{\mathcal C}(x,y). \]

The following definition of weak equivalence of simplicial
categories is a natural generalization of the notion of
equivalence of categories.

\begin{definition} \label{dkequiv}
A simplicial functor $f:\mathcal C \rightarrow \mathcal D$ is a
\emph{Dwyer-Kan equivalence} or \emph{DK-equivalence} if the
following two conditions hold:
\begin{enumerate}
\item For any objects $x$ and $y$ of $\mathcal C$, the induced map
$\Map(x,y) \rightarrow \Map(fx, fy)$ is a weak equivalence of
simplicial sets.

\item The induced map on the categories of components $\pi_0 f
\colon \pi_0 \mathcal C \rightarrow \pi_0 \mathcal D$ is an
equivalence of categories.
\end{enumerate}
\end{definition}

The idea of obtaining a simplicial category from a model category
$\mathcal M$ goes back to several papers of Dwyer and Kan
\cite{dkcalc}, \cite{dkfncxes}, \cite{dksimploc}.  In fact, they
define two different methods of doing so, the simplicial
localization $L\mathcal M$ \cite{dksimploc} and the hammock
localization $L^H \mathcal M$ \cite{dkfncxes}.  The first has the
advantage of being easier to describe, while the second is more
convenient for making calculations.

It should be noted that these constructions can be made for more
general categories with weak equivalences, and do not depend on
the model structure if we are willing to ignore the potential
set-theoretic difficulties.  However, as with the homotopy
category construction, the hammock localization in particular can
be defined much more nicely when we have the additional structure
of a model category.

We begin with the construction of the simplicial localization
$L\mathcal M$.  Recall that, given a category $\mathcal M$ with
some choice of weak equivalences $\mathcal W$, we denote the
localization of $\mathcal M$ with respect to $\mathcal W$ by
$\mathcal W^{-1} \mathcal M$.  This localization is obtained from
$\mathcal M$ by formally inverting the maps of $\mathcal W$.
Further, recall that, given a category $\mathcal M$, we denote by
$F \mathcal M$ the \emph{free category} on $\mathcal M$, or
category with the same objects as $\mathcal M$ and morphisms
freely generated by the non-identity morphisms of $\mathcal M$.
Note in particular that there are natural functors $F \mathcal M
\rightarrow \mathcal M$ and $F \mathcal M \rightarrow F^2 \mathcal
M$ \cite[2.4]{dksimploc}.  These functors can be used to define a
simplicial resolution $F_*\mathcal M$, which is a simplicial
category with the category $F^{k+1} \mathcal M$ at level $k$
\cite[2.5]{dksimploc}.

We can apply this same construction to the subcategory $\mathcal
W$ to obtain a simplicial resolution $F_* \mathcal W$.  Using
these two resolutions, we have the following definition.

\begin{definition} \cite[4.1]{dksimploc}
The \emph{simplicial localization} of $\mathcal M$ with respect to
$\mathcal W$ is the localization $(F_* \mathcal W)^{-1}
(F_*\mathcal M)$.  This simplicial localization is denoted
$L(\mathcal M, \mathcal W)$ or simply $L \mathcal M$.
\end{definition}

The following result gives interesting information about the
mapping spaces in $L\mathcal M$ in the case where $\mathcal W$ is
all of $\mathcal M$.

\begin{prop} \cite[5.5]{dksimploc} \label{groupoid}
Suppose that $\mathcal W = \mathcal M$ and $\nerve(\mathcal M)$ is
connected.
\begin{enumroman}
\item The simplicial localization $L \mathcal M$ is a simplicial
groupoid, so for all objects $x$ and $y$, the simplicial sets
$\Map_{L \mathcal M}(x,y)$ are all isomorphic. In particular, the
simplicial sets $\Map_{L \mathcal M}(x,x)$ are all isomorphic
simplicial groups.

\item The classifying complex $B \Map_{L \mathcal M}(x,x)$ has the
homotopy type of $\nerve(\mathcal M)$, and thus each simplicial
set $\Map_{L \mathcal M}(x,y)$ has the homotopy type of the loop
space $\Omega (\nerve(\mathcal M))$.
\end{enumroman}
\end{prop}

We now turn to the other construction, that of the hammock
localization.  Again, let $\mathcal M$ be a category with a
specified subcategory $\mathcal W$ of weak equivalences.

\begin{definition} \cite[3.1]{dkfncxes}
The \emph{hammock localization} of $\mathcal M$ with respect to
$\mathcal W$, denoted $L^H(\mathcal M, \mathcal W)$, or simply
$L^H\mathcal M$, is the simplicial category defined as follows:
\begin{enumerate}
\item The simplicial category $L^H \mathcal M$ has the same
objects as $\mathcal M$.

\item Given objects $X$ and $Y$ of $\mathcal M$, the simplicial
set $\Map_{L^H \mathcal M}(X,Y)$ has as $k$-simplices the reduced
hammocks of width $k$ and any length between $X$ and $Y$, or
commutative diagrams of the form
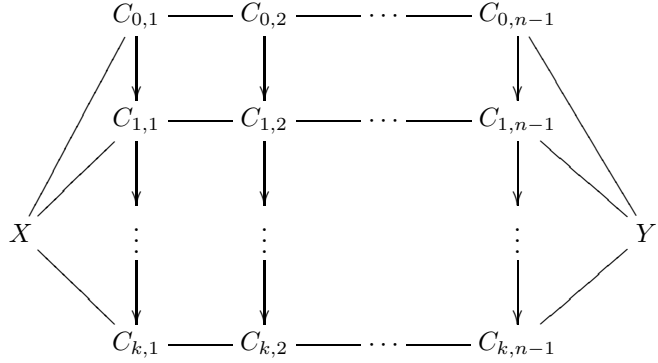
\[ \xymatrix{& C_{0,1} \ar[d] \ar@{-}[r] & C_{0,2} \ar[d] \ar@{-}[r] & \cdots \ar@{-}[r] &
C_{0,n-1} \ar[d] & \\
& C_{1,1} \ar[d] \ar@{-}[r] & C_{1,2} \ar[d] \ar@{-}[r] & \cdots
\ar@{-}[r] & C_{1,n-1}
\ar[d] & \\
X \ar@{-}[uur] \ar@{-}[ur] \ar@{-}[dr] & \vdots \ar[d] & \vdots
\ar[d] && \vdots \ar[d] &
Y \ar@{-}[uul] \ar@{-}[ul], \ar@{-}[dl] \\
& C_{k,1} \ar@{-}[r] & C_{k,2} \ar@{-}[r] & \cdots \ar@{-}[r] &
C_{k,n-1} & }
\] in which
\begin{enumroman}
\item the length of the hammock is any integer $n \geq 0$,

\item the vertical maps are all in $\mathcal W$,

\item in each column all the horizontal maps go the same
direction, and if they go to the left, then they are in $\mathcal
W$,

\item the maps in adjacent columns go in opposite directions, and

\item no column contains only identity maps.
\end{enumroman}
\end{enumerate}
\end{definition}

\begin{prop} \cite[2.2]{dkcalc}
For a given model category $\mathcal M$, the simplicial categories
$L \mathcal M$ and $L^H \mathcal M$ are DK-equivalent.
\end{prop}

We should add that the description of the hammock localization can
be greatly simplified if we make use of the model category
structure on $\mathcal M$.  In this case, Dwyer and Kan prove that
it suffices to consider hammocks of length 3 such as the following
\cite[\S 8]{dkcalc}:
\[ \xymatrix@1{X & C_{0,1} \ar[l]_\simeq \ar[r] &
C_{0,2}  &  Y \ar[l]_-\simeq}. \]

Restricting to the category of simplicial categories with a fixed
set $\mathcal O$ of objects, Dwyer and Kan prove the existence of
a model structure on this category, which we denote
$\mathcal{SC}_\mathcal O$ \cite[7.2]{dksimploc}.  In this
situation, the weak equivalences are the DK-equivalences, but with
the objects fixed the second condition follows immediately from
the first. The fibrations in this model structure are given by the
functors $f\colon \mathcal C \rightarrow \mathcal D$ inducing, for
any objects $x$ and $y$, fibrations of simplicial sets
$\Map_\mathcal C(x,y) \rightarrow \Map_\mathcal D(x,y)$.

The cofibrations are then defined to be the maps with the left
lifting property with respect to the acyclic fibrations.  However,
they can be more precisely characterized.  To do so, we recall the
definition of a free map of simplicial categories.

\begin{definition}\cite[7.4]{dksimploc}
A map $f:\mathcal C \rightarrow \mathcal D$ in $\sco$ is
\emph{free} if
\begin{enumerate}
\item $f$ is a monomorphism,

\item if $\ast$ denotes the free product, then in each simplicial
dimension $k$, the category $\mathcal D_k$ admits a unique free
factorization $\mathcal D_k= f(\mathcal C_k) \ast \mathcal F_k$,
where $\mathcal F_k$ is a free category, and

\item for each $k \geq 0$, all degeneracies of generators of
$\mathcal F_k$ are generators of $\mathcal F_{k+1}$.
\end{enumerate}
\end{definition}

\begin{definition}\cite[7.5]{dksimploc}
A map $f:\mathcal C \rightarrow \mathcal D$ of simplicial
categories is a \emph{strong retract} of a map $f':\mathcal C
\rightarrow \mathcal D'$ if there exists a commutative diagram
\[ \xymatrix{& \mathcal C \ar[ldd]_f \ar[rdd]^f \ar[d]^{f'}& \\
& \mathcal D' \ar[dr] & \\
\mathcal D \ar[ur] \ar[rr]^{id} && \mathcal D} \]
\end{definition}

Using these definitions, Dwyer and Kan prove the following result.

\begin{prop} \cite[7.6]{dksimploc}
The cofibrations of $\sco$ are precisely the strong retracts of
free maps. In particular, a cofibrant simplicial category is a
retract of a free category.
\end{prop}

This result can then be generalized to the category of all
simplicial categories, in which the DK-equivalences are the weak
equivalences.  If $\mathcal C$ is a simplicial category, a
morphism $e \in \Hom_\mathcal C(a,b)_0$ is a \emph{homotopy
equivalence} if it becomes an isomorphism in $\pi_0 \mathcal C$.

\begin{theorem} \cite[1.1]{simpcat}
There is a model category structure $\mathcal{SC}$ on the category
of small simplicial categories in which
\begin{itemize}
\item the weak equivalences are the Dwyer-Kan equivalences, and

\item the fibrations are the maps $f: \mathcal C \rightarrow
\mathcal D$ satisfying the following two conditions:
\begin{enumroman}
\item For any objects $a_1$ and $a_2$ in $\mathcal C$, the map
\[ \Hom_\mathcal C (a_1,a_2) \rightarrow \Hom_\mathcal D (fa_1,fa_2) \]
is a fibration of simplicial sets.

\item For any object $a_1$ in $\mathcal C$, $b$ in $\mathcal D$,
and homotopy equivalence $e:fa_1 \rightarrow b$ in $\mathcal D$,
there is an object $a_2$ in $\mathcal C$ and homotopy equivalence
$d:a_1 \rightarrow a_2$ in $\mathcal C$ such that $fd=e$.
\end{enumroman}
\end{itemize}
\end{theorem}

\section{Complete Segal Spaces}

Here we define complete Segal spaces and describe Rezk's model
structure on the category of simplicial spaces, in which the
complete Segal spaces are the fibrant-cofibrant objects.

Recall that by a simplicial space we mean a simplicial object in
the category of simplicial sets, or functor $\Deltaop \rightarrow
\SSets$.  In section 2, we described the Reedy model category
structure on this category, in which both the weak equivalences
and cofibrations are defined levelwise.  The model structure
$\css$ is given by a localization of this structure with respect
to a set of maps.

We begin with the definition of a Segal space. In
\cite[4.1]{rezk}, Rezk defines for each $0 \leq i \leq n-1$ a map
$\alpha^i \colon [1] \rightarrow [n]$ in ${\bf \Delta}$ such that
$0 \mapsto i$ and $1 \mapsto i+1$. There is a corresponding map
$\alpha_i \colon \Delta[1] \rightarrow \Delta[n]$.  Then for each
$n$ he defines the simplicial space
\[ G(n)^t= \bigcup_{i=0}^{n-1} \alpha_i \Delta [1]^t \subset \Delta
[n]^t. \]

Let $X$ be a Reedy fibrant simplicial space.  There is a weak
equivalence of simplicial sets
\[ \Map_{\SSets_c^{\Deltaop}} (G(n)^t, X) \rightarrow \underbrace{X_1 \times_{X_0} \cdots
\times_{X_0} X_1}_n, \] where the right hand side is the limit of
the diagram
\[ \xymatrix{X_1 \ar[r]^{d_0} & X_0 & X_1 \ar[l]_{d_1}
\ar[r]^{d_0} & \ldots \ar[r]^{d_0} & X_0 & X_1 \ar[l]_{d_1}} \]
with $n$ copies of $X_1$.

Now, given any $n$, define the map $\varphi^n \colon G(n)^t
\rightarrow \Delta [n]^t$ to be the inclusion map.  Then for any
Reedy fibrant simplicial space $W$ there is a map
\[ \varphi_n = \Map_{\SSets_c^{\Deltaop}}(\varphi^n, W) \colon \Map_{\SSets_c^{\Deltaop}}(\Delta
[n]^t,W) \rightarrow \Map_{\SSets_c^{\Deltaop}}(G(n)^t,W). \] More
simply written, this map is
\[ \varphi_n \colon W_n \rightarrow \underbrace{W_1 \times_{W_0} \cdots
\times_{W_0} W_1}_n \] and is often called a \emph{Segal map}. The
Segal map is actually defined for any simplicial space $W$, but
here we assume Reedy fibrancy so that the mapping spaces involved
are homotopy invariant.

\begin{definition} \cite[4.1]{rezk}
A Reedy fibrant simplicial space $W$ is a \emph{Segal space} if
for each $n \geq 2$ the Segal map
\[ \varphi_n: W_n \rightarrow W_1 \times_{W_0} \cdots \times_{W_0}
W_1 \] is a weak equivalence of simplicial sets.
\end{definition}

In fact, there is a model category structure $\sesp$ on the
category of simplicial spaces in which the fibrant objects are
precisely the Segal spaces \cite[7.1]{rezk}.  This model structure
is obtained from the Reedy structure via localization.

The idea is that in a Segal space there is a notion of
``composition," at least up to homotopy.  In fact, given a Segal
space, we can sensibly use many categorical notions.  We summarize
some of these ideas here; a detailed description is given by Rezk
\cite{rezk}.  The \emph{objects} of a Segal space $W$ are given by
the set $W_{0,0}$.  Given the map
\[ (d_1, d_0): W_1 \rightarrow W_0 \times W_0, \] the \emph{mapping
space} $\map_W(x,y)$ is given by the fiber of this map over
$(x,y)$.  (The fact that $W$ is Reedy fibrant guarantees that this
mapping space is homotopy invariant.)  Two maps $f,g \in
\map_W(x,y)_0$ are \emph{homotopic} if they lie in the same
component of the simplicial set $\map_W(x,y)$.  Thus, we define
the space of homotopy equivalences $W_{\hoequiv} \subseteq W_1$ to
consist of all the components containing homotopy equivalences.

Given any $(x_0, \ldots ,x_n) \in W_{0,0}^{n+1}$, let $\map_W(x_0,
\ldots ,x_n)$ denote the fiber of the map
\[ (\alpha_0, \ldots ,\alpha_n) \colon W_n \rightarrow W_0^{n+1}
\] over $(x_0, \ldots, x_n)$.  Consider the commutative diagram
\[ \xymatrix{W_n = \Map(\Delta [n]^t, W) \ar[rr]^{\varphi_k} \ar[dr]
&& \Map(G(n)^t,W) \ar[dl] \\
& W_0^{n+1} & } \] and notice that, since $W$ is a Segal space,
the horizontal arrow is a weak equivalence and a fibration.  In
particular, this map induces an acyclic fibration on the fibers of
the two vertical arrows,
\[ \map_W(x_0, \ldots, x_n) \rightarrow \map_W(x_{n-1},x_n) \times
\cdots \times \map_W(x_0, x_1). \]

Given $f \in \map(x,y)_0$ and $g \in \map(y,z)_0$, their
\emph{composite} is a lift of $(g,f) \in \map(y,z) \times
\map(x,y)$ along $\varphi_2$ to some $k \in \map(x,y,z)_0$.  The
\emph{result} of this composition is defined to be $d_1(k) \in
\map(x,z)_0$.  It can be shown that any two results are homotopic,
so we can use $g \circ f$ unambiguously.

Then, the \emph{homotopy category} of $W$, denoted $\Ho(W)$, has
as objects the elements of the set $W_{0,0}$, and
\[ \Hom_{\Ho(W)}(x,y) = \pi_0 \map_W(x,y). \]  A \emph{homotopy
equivalence} in $W$ is a 0-simplex of $W_1$ whose image in
$\Ho(W)$ is an isomorphism.

\begin{definition}
A map $f \colon W \rightarrow Z$ of Segal spaces is a
\emph{Dwyer-Kan equivalence} if
\begin{enumerate}
\item for any objects $x$ and $y$ of $W$, the induced map
$\map_W(x,y) \rightarrow \map_Z(fx,fy)$ is a weak equivalence of
simplicial sets, and

\item the induced map $\Ho(W) \rightarrow \Ho(Z)$ is an
equivalence of categories.
\end{enumerate}
\end{definition}

Notice that the definition of these maps bears a striking
resemblance to that of the Dwyer-Kan equivalences between
simplicial categories, hence the use of the same name.

For a Segal space $W$, note that the degeneracy map $s_0 \colon
W_0 \rightarrow W_1$ factors through the space of homotopy
equivalences $W_{\hoequiv}$, since the image of $s_0$ consists of
``identity maps."  Given this fact, we are now able to give a
definition of complete Segal space.

\begin{definition} \cite[\S 6]{rezk} \label{css}
A Segal space $W$ is a \emph{complete Segal space} if the map $W_0
\rightarrow W_{\hoequiv}$ given above is a weak equivalence of
simplicial sets.
\end{definition}

The idea behind this notion is that, although $W_0$ is not
required to be discrete, as the objects are for a simplicial
category, it is not heuristically too different from a simplicial
space with discrete 0-space.  (This viewpoint is further confirmed
by the comparison of complete Segal spaces with Segal categories,
which are essentially the analogues of Segal spaces with discrete
0-space \cite[6.3]{thesis}.)

Now, we give a description of the model structure $\css$.  We do
not give all the details here, such as a description of an
arbitrary weak equivalence, but refer the interested reader to
Rezk's paper \cite[\S 7]{rezk}.

\begin{theorem} \cite[7.2]{rezk}
There is a model structure $\css$ on the category of simplicial
spaces, obtained as localization of the Reedy model structure,
such that
\begin{enumerate}
\item the fibrant objects are precisely the complete Segal spaces,

\item the cofibrations are the monomorphisms; in particular, every
object is cofibrant,

\item the weak equivalences between Segal spaces are Dwyer-Kan
equivalences, and

\item the weak equivalences between complete Segal spaces are
levelwise weak equivalences of simplicial sets.
\end{enumerate}
Furthermore, $\css$ has the additional structure of a simplicial
model category and is cartesian closed.
\end{theorem}

The fact that $\css$ is cartesian closed allows us to consider,
for any complete Segal space $W$ and simplicial space $X$, the
complete Segal space $W^X$.  In particular, using the simplicial
structure, the simplicial set at level $n$ is given by
\[ (W^X)_n = \Map(X \times \Delta[n]^t, W). \]  If $W$ is a (not
necessarily complete) Segal space, then $W^X$ is again a Segal
space; in other words, the model category $\sesp$ is also
cartesian closed.

We denote the functorial fibrant replacement functor in $\css$ by
$L_{\css}$. Thus, given any simplicial space $X$, there is a
weakly equivalent complete Segal space $L_{\css}X$.

This model structure is connected to the model structure
$\mathcal{SC}$ by a chain of Quillen equivalences as follows. Each
of these model categories is Quillen equivalent to a model
structure $\Secat_f$ on the category of Segal precategories. A
\emph{Segal precategory} is a simplicial space $X$ with $X_0$ a
discrete space.  A \emph{Segal category} is then a Segal
precategory with the Segal maps weak equivalences.  In the model
structure $\Secat_f$, the fibrant objects are Segal categories,
and so it is considered a Segal category model structure on the
category of Segal precategories. We have the following chain of
Quillen equivalences, with the left adjoint functors topmost:
\[ \mathcal{SC} \leftrightarrows \Secat_f \rightleftarrows \css.
\]
The right adjoint $\mathcal{SC} \rightarrow \Secat_f$ is given by
the nerve functor, and the left adjoint $\Secat_f \rightarrow
\css$ is given by the inclusion functor.

There is actually another chain of Quillen equivalences connecting
the two model structures; in this case, both $\mathcal{SC}$ and
$\css$ are Quillen equivalent to Joyal's model structure
$\mathcal{QC}at$ on the category of simplicial sets
\cite{joyalquasi}. The fibrant objects in $\mathcal{QC}at$ are
\emph{quasi-categories}, or simplicial sets $K$ such that a dotted
arrow lift exists making the diagram
\[ \xymatrix{V[m,k] \ar[r] \ar[d] & K \\
\Delta[m] \ar@{-->}[ur] & } \] commute for any $0 < k < m$. The
chain of Quillen equivalences in this case is given by
\[ \mathcal{SC} \rightleftarrows \mathcal{QC}at \leftrightarrows
\css. \]  The right adjoint $\mathcal{SC} \rightarrow
\mathcal{QC}at$ is given by Cordier and Porter's coherent nerve
functor \cite{cp}, \cite[2.10]{joyalsc}, and the right adjoint
$\css \rightarrow \mathcal{QC}at$ is given given by sending a
simplicial space $W$ to the simplicial set $W_{\ast, 0}$
\cite[4.11]{jt}.  It is a consequence of work of Joyal \cite[\S
1-2]{joyalsc} and of Joyal and Tierney \cite[\S 4-5]{jt} that the
simplicial space obtained from a simplicial category via these
functors is weakly equivalent to the one obtained from the
composite functor described in the previous paragraph.

\section{Obtaining Complete Segal Spaces from Simplicial
Categories and Model Categories}

In this section, we describe several different ways of obtaining a
complete Segal space.  First, we look at a particularly nice
functor which Rezk uses to modify the notion of a nerve of a
category. Then we look at how this functor can be generalized to
one on any simplicial category, and how a similar idea can be used
to get a complete Segal space from any model category.  We then
consider the functors used in the Quillen equivalences connecting
$\mathcal{SC}$ and $\mathcal{CSS}$.

Let us begin with Rezk's classifying diagram construction, which
associates to any small category $\mathcal C$ a complete Segal
space $N\mathcal C$.  First, we denote by $\nerve (\mathcal C)$
the ordinary nerve, which is the simplicial set given by
$(\nerve(\mathcal C))_n= \Hom([n], \mathcal C)$. Further, we
denote by $\iso (\mathcal C)$ the maximal subgroupoid of $\mathcal
C$, or subcategory of $\mathcal C$ with all objects of $\mathcal
C$ and whose only morphisms are the isomorphisms of $\mathcal C$.
By $\mathcal C^{[n]}$ we denote the category of functors $[n]
\rightarrow \mathcal C$, or the category whose objects are
$n$-chains of composable morphisms in $\mathcal C$.

\begin{definition} \cite[3.5]{rezk}
The \emph{classifying diagram} $N \mathcal C$ is the simplicial
space given by $(N \mathcal C)_n = \nerve(\iso (\mathcal
C^{[n]}))$.
\end{definition}

Thus, $(N \mathcal C)_0$ is simply the nerve of $\iso (\mathcal
C)$, and $(N \mathcal C)_1$ is the nerve of the maximal
subgroupoid of the morphism category of $\mathcal C$.  In
particular, information about invertible morphisms of $\mathcal C$
is encoded at level 0, while information about the other morphisms
of $\mathcal C$ does not appear until level 1.

Thus, the classifying diagram of a category can be regarded as a
more refined version of the nerve, since, unlike the ordinary
nerve construction, it enables one to recover information about
whether morphisms are invertible or not.  This construction is
also particularly useful for our purposes due to the following
result.

\begin{prop} \cite[6.1]{rezk}
If $\mathcal C$ is a small category, then its classifying diagram
$N \mathcal C$ is a complete Segal space.
\end{prop}

However, this construction, as defined above, cannot be used to
assign a complete Segal space to any simplicial category, since,
beginning with level 1, we would have homotopy invariance problems
with a simplicial set of objects in $\mathcal C^{[1]}$.  Rezk
defines an analogous functor, though, from the category of small
simplicial categories which is similar in spirit to the
classifying diagram but avoids these difficulties.

Let $I[m]$ denote the category with $m+1$ objects and a single
isomorphism between any two objects, and let $E(m)=
\nerve(I[m])^t$.  If $W$ is a Segal space and $X$ is any
simplicial space, recall that $W^X$ denotes the internal hom
object, which is a Segal space.  With these notations in place, we
can give the definition of Rezk's completion functor. Let $W$ be a
Segal space.  Then its completion $\widehat W$ is defined as a
fibrant replacement in $\css$ of the simplicial space $\widetilde
W$ defined by
\[ \widetilde W_n= \diag ([m] \mapsto
\Map_{SSets^{\Deltaop}}(E(m), W^{\Delta[n]^t})) = \diag ([m]
\mapsto (W^{E(m)})_n). \]

From a simplicial category $\mathcal C$, then, we can takes its
nerve to obtain a simplicial space, followed by a fibrant
replacement functor in the Segal space model structure, to obtain
a Segal space $W$.  From $W$ we can then pass to a complete Segal
space via this completion functor.  We will denote this complete
Segal space $L_C(W)$, or $L_C(\mathcal C)$ where $W$ comes from a
simplicial category as just described.

The first important fact about this completion functor is that the
completion map $i_W:W \rightarrow \widehat W = L_C(W)$ is not only
a weak equivalence in the model category $\css$, but is also a
Dwyer-Kan equivalence of Segal spaces \cite[\S 14]{rezk}.
Furthermore, this completion functor restricts nicely to the
classifying diagram in the case where $\mathcal C$ is a discrete
category.

\begin{prop} \cite[14.2]{rezk}
If $\mathcal C$ is a discrete category, then $L_C (\mathcal C)$ is
isomorphic to $N \mathcal C$.
\end{prop}

If we begin with a model category $\mathcal M$ with subcategory of
weak equivalences $\mathcal W$, a functor analogous to the
classifying diagram functor can be used to obtain a complete Segal
space.  In this case, rather than taking the subcategory $\iso
(\mathcal M)$ of isomorphisms of $\mathcal M$, we take the
subcategory of weak equivalences, denoted $\we(\mathcal M)$. Thus,
Rezk defines the \emph{classification diagram} of $(\mathcal M,
\mathcal W)$, denoted $N(\mathcal M, \mathcal W)$, by
\[ N(\mathcal M, \mathcal W)_n= \nerve (\we (\mathcal M^{[n]})). \]
Unlike the classifying diagram, the classification diagram of a
model category is not necessarily a complete Segal space as
stated, but taking a Reedy fibrant replacement of it results in a
complete Segal space, as we show in the next section.

Lastly, we have the two functors given by the two different chains
of Quillen equivalences between the model categories
$\mathcal{SC}$ and $\css$.  As mentioned in the previous section,
these two functors are equivalent.  In each case, the resulting
simplicial space is not Reedy fibrant in general, and so not a
complete Segal space, but applying the fibrant replacement functor
$L_{\css}$ results in a complete Segal space.

The first of these composite functors, in particular, is simple to
describe abstractly, as in the previous section, but it has a
disadvantage over Rezk's functor in that it gives very little
insight into what the resulting complete Segal space looks like.
In the next section, we prove that the two functors from
$\mathcal{SC}$ to $\css$ result in weakly equivalent complete
Segal spaces, and that if we use Rezk's classification diagram
construction we get a weakly equivalent complete Segal space to
the one we would obtain by taking the simplicial localization
followed by his completion functor.  We then use Rezk's functor to
describe what the complete Segal space corresponding to a
simplicial category looks like.

\section{Comparison of Functors from $\mathcal{SC}$ to $\mathcal{CSS}$}

Here we prove that each of the functors we have described all give
rise to complete Segal spaces weakly equivalent to those given in
the previous section.  We begin by stating the result that
establishes the equivalence between Rezk's completion functor $L_C
\colon \mathcal{SC} \rightarrow \css$ and the functor arising from
the chain of Quillen equivalences factoring through $\Secat_f$.
Let $L_{\css}$ denote the functorial fibrant replacement functor
in $\css$.

\begin{theorem}
If $\mathcal C$ is a simplicial category, then the complete Segal
spaces $L_C(\mathcal C)$ and $L_{\css}(\nerve(\mathcal C))$ are
weakly equivalent in $\css$.
\end{theorem}

\begin{proof}
Let $L_S$ denote a fibrant replacement functor in the Segal space
model structure $\sesp$ on the category of simplicial spaces.  The
fact that the two functors in question result in weakly equivalent
complete Segal spaces can be shown by considering the following
chain of weak equivalences:
\[ L_{\css}(\nerve(\mathcal C)) \leftarrow \nerve(\mathcal C)
\rightarrow L_S \nerve(\mathcal C) \rightarrow L_C(\mathcal C).
\]  The map on the left is the localization functor in $\css$ and
so is a weak equivalence in $\css$.  The middle map is a weak
equivalence in $\sesp$ and therefore also a weak equivalence in
$\css$, since the latter model category is a localization of the
former. The map on the right is Rezk's completion, and it is a
weak equivalence in $\css$, as given in the previous section.
Therefore, the objects at the far left and right of this zigzag,
both of which are complete Segal spaces, are weakly equivalent as
objects of $\css$.
\end{proof}

Now, we would like to compare either of these functors to the
classifying diagram construction for a model category $\mathcal
M$.  In other words, we want to show that $N(\mathcal M, \mathcal
W)$ is equivalent to $L_C(L^H \mathcal M)$, where we first take
the hammock localization of $\mathcal M$ to obtain a simplicial
category, and then apply Rezk's functor $L_C$.

An initial problem here is that $N(\mathcal M, \mathcal W)$ is not
necessarily Reedy fibrant, and so it is not necessarily a complete
Segal space.  We prove that a Reedy fibrant replacement of it,
denoted $N(\mathcal M, \mathcal W)^f$, is in fact a complete Segal
space in the process of comparing the ``mapping spaces" in this
Reedy fibrant replacement to the mapping spaces of the hammock
localization $L^H\mathcal M$.

\begin{theorem} \label{equiv}
Let $\mathcal M$ be a model category, and let $\mathcal W$ denote
its subcategory of weak equivalences.  Then $N(\mathcal M,
\mathcal W)^f$, is a complete Segal space. Furthermore, for any
objects $x,y$ of $\mathcal M$, there is a weak equivalence of
spaces $map_{N(\mathcal M, \mathcal W)^f}(x,y) \simeq
\Map_{L^H\mathcal M}(x,y)$, and there is an equivalence of
categories $\Ho(N(\mathcal M, \mathcal W)^f) \approx \Ho(\mathcal
M)$.
\end{theorem}

This result was proved by Rezk in the case where $\mathcal M$ is a
simplicial model category \cite[8.3]{rezk}, namely, in the case
where we do not need to pass to the simplicial localization of
$\mathcal M$ to consider its function complexes.  However, here we
prove that, as he conjectured \cite[8.4]{rezk}, the result
continues to hold in this more general case.

We prove this theorem very similarly to the way Rezk proves it in
the more restricted case, using a proposition of Dwyer and Kan. To
begin, we introduce some terminology.  Let $\mathcal M$ be a model
category. A \emph{classification complex} of $\mathcal M$, as
defined in \cite[1.2]{dkclass}, is the nerve of any subcategory
$\mathcal C$ of $\mathcal M$ such that
\begin{enumerate}
\item every map in $\mathcal C$ is a weak equivalence,

\item if $f \colon X \rightarrow Y$ in $\mathcal M$ is a weak
equivalence and either $X$ or $Y$ is in $\mathcal C$, then $f$ is
in $\mathcal C$, and

\item $\nerve(\mathcal C)$ is homotopically small; i.e., each
homotopy group of $|\nerve(\mathcal C)|$ is small
\cite[2.2]{dkfncxes}.
\end{enumerate}

The \emph{special classification complex} $sc(X)$ of an object $X$
in $\mathcal M$ is a connected classification complex containing
$X$.

Let $\mathcal M$ be a model category and $X$ a fibrant-cofibrant
object of $\mathcal M$.  Denote by $Aut^h(X)$ the simplicial
monoid of weak equivalences given by $Aut^h_{L^H \mathcal M}(X)$
in the hammock localization $L^H \mathcal M$, and by $B \Aut^h(X)$
its classifying complex.

The following proposition was proved by Dwyer and Kan in
\cite[2.3]{dkclass} in the case that $\mathcal M$ is a simplicial
model category. However, the proof does not actually require the
simplicial structure; in fact, their proof is essentially the one
given below, with the extra step showing that the mapping spaces
in the hammock localization are equivalent to those given by the
simplicial structure of $\mathcal M$ \cite[4.8]{dkfncxes}.

\begin{prop} \label{sc}
Let $X$ be an object of a model category $\mathcal M$.  The
classifying complex $B Aut^h(X)$ is weakly equivalent to the
special classification complex of $X$, $\text{sc}(X)$, and the two
can be connected by a finite zig-zag of weak equivalences.
\end{prop}

\begin{proof}
Let $\mathcal W$ be the subcategory of weak equivalences of
$\mathcal M$.  Consider the connected component of
$\nerve(\mathcal W)$ containing $X$.  For the rest of this proof,
we assume that $\mathcal W$ is such that its nerve is connected.
We further assume that $\nerve(\mathcal W)$ is homotopically
small, taking an appropriate subcategory, as described in
\cite[2.3]{dkfncxes}, if necessary.

In this case, by Proposition \ref{groupoid}, the function
complexes $\Map_{L \mathcal W}(X,X)$ are all isomorphic.
Furthermore, by the same result, the classifying complex $B
\Map_{L \mathcal W}(X,X)$ has the homotopy type of
$\nerve(\mathcal W)$. Thus, we can take $\nerve(\mathcal W)$ as
$sc(X)$.

Now, as in the statement of the proposition, we take $\Aut^h(X)$
to consist of the components of $\Map_{L^H \mathcal M}(X,X)$ which
are invertible in $\pi_0 \map{L^H \mathcal M}(X,X)$.  But, by
\cite[4.6(ii)]{dkfncxes}, the map $B \Map_{L^H \mathcal W}(X,X)
\rightarrow B \Aut^h(X)$ is a weak equivalence of simplicial sets.
Since $L^H \mathcal W$ can be connected to $L \mathcal W$ by a
finite string of weak equivalences, it follows that so can
$\Map_{L^H \mathcal W}(X,X)$ and $\Map_{L \mathcal W} (X,X)$.
Thus, $B \Map_{L \mathcal W}(X,X)$ and $B \Aut^h(X)$ can also be
connected by such a string. It follows that $sc(X)$ has the same
homotopy type as $B \Aut^h(X)$.
\end{proof}

\begin{proof}[Proof of Theorem \ref{equiv}]
Consider the category $\mathcal M^{[n]}$ of functors $[n]
\rightarrow \mathcal M$.  If $\mathcal M$ is a model category,
then $\mathcal M^{[n]}$ can be given the structure of a model
category with the weak equivalences and fibrations given by
levelwise weak equivalences and fibrations in $\mathcal M$.  Given
any map $[m] \rightarrow [n]$, we obtain a functor $\mathcal
M^{[n]} \rightarrow \mathcal M^{[m]}$.

Let $Y= (y_0 \rightarrow y_1 \cdots \rightarrow y_n)$ be a
fibrant-cofibrant object of $\mathcal M^{[n]}$.  It restricts to
an object $Y'=(y_0 \rightarrow y_1 \cdots \rightarrow y_{n-1})$ in
$\mathcal M^{[n-1]}$.  From this map, we obtain a map of
simplicial sets
\[ B \Aut^h_{L^H\mathcal M^{[n]}}(Y) \rightarrow B
\Aut^h_{L^H \mathcal M} (y_n) \times B \Aut^h_{L^H \mathcal
M^{[n-1]}}(Y').
\] The homotopy fiber of this map is weakly equivalent to the
union of those components of $\Map_{L^H \mathcal M}(y_{n-1}, y_n)$
containing the conjugates of the map $f_{n-1}:y_{n-1} \rightarrow
y_n$, or maps $j \circ f_{n-1} \circ i$, where $i$ and $j$ are
self-homotopy equivalences.

Iterating this process, we can take the homotopy fiber of the map
\[ B \Aut^h_{L^H \mathcal M^{[n]}}(Y) \rightarrow B \Aut^h_{L^H
\mathcal M} (y_n) \times \cdots \times B \Aut^h_{L^H \mathcal M}
(y_0), \] which is weakly equivalent to the union of the
components of
\[ \Map_{L^H \mathcal M}(y_{n-1},y_n) \times \cdots \times \Map_{L^H \mathcal
M}(y_0,y_1) \] containing conjugates of the sequence of maps
$f_i:y_i \rightarrow y_{i+1}$, $0 \leq i \leq n-1$.  However,
applying Proposition \ref{sc} to the map in question shows that
this simplicial set is also the homotopy fiber of the map
\[ \text{sc}(Y) \rightarrow \text{sc}(y_n) \times \cdots \times
\text{sc}(y_0). \]

Let $U$ denote the simplicial space $N(\mathcal M, \mathcal W)$ so
that $U_n= \nerve(\we(\mathcal M^{[n]}))$.  Then, let $V$ be a
Reedy fibrant replacement of $U$, from which we get weak
equivalences $U_n \rightarrow V_n$ for all $n \geq 0$.

For each $n \geq 0$, there exists a map $p_n:U_n \rightarrow
U_0^{n+1}$ given by iterated face maps to the ``objects."  Then,
for every $(n+1)$-tuple of objects $(x_0, x_1, \ldots, x_n)$, the
homotopy fiber of $p_n$ over $(x_0, \ldots ,x_n)$, given by
\[ \map_V(x_{n-1},x_n) \times \cdots \times \map_V(x_0, x_1), \]
is weakly equivalent to
\[ \Map_{L^H \mathcal M} (x_{n-1}^{cf},x^{cf}_n) \times \cdots \times \Map_{L^H \mathcal
M}(x^{cf}_0, x^{cf}_1) \] where $x^{cf}$ denotes a
fibrant-cofibrant replacement of $X$ in $\mathcal M$.  It follows
that once we take the Reedy fibrant replacement $V$ of $U$, it is
a Segal space.

Now, consider the set $\pi_0 U_0$, which consists of the weak
equivalence classes of objects in $\mathcal M$; it follows that
$\pi_0 V_0$ is an isomorphic set.  Further, note that
\[ \Hom_{\Ho(\mathcal M)}(x,y)= \pi_0 \Map_{L^H \mathcal
M}(x^{cf}, y^{cf}). \]  Thus, we have shown that $\Ho(\mathcal M)$
is equivalent to $\Ho(V)$.

It remains to show that $V$ is a complete Segal space. Consider
the space $V_\hoequiv \subseteq V_1$, and define $U_\hoequiv$ to
be the preimage of $V_\hoequiv$ under the natural map $U
\rightarrow V$.  Since $V$ is a Reedy fibrant replacement for $U$,
it suffices to show that the complete Segal space condition holds,
i.e., that $U_0 \rightarrow U_\hoequiv$ is a weak equivalence of
simplicial sets.  Notice that $U_\hoequiv$ must consist precisely
of the components of $U_1$ whose 0-simplices come from weak
equivalences in $\mathcal M$.  In other words, $U_\hoequiv =
\nerve(\we(\we(\mathcal M))^{[1]})$.

There is an adjoint pair of functors
\[ F \colon \mathcal M^{[1]} \rightleftarrows \mathcal M \colon G
\] for which $F(x \rightarrow y)= x$ and $G(x)= \text{id}_x$.
This adjoint pair can be restricted to an adjoint pair
\[ F \colon \nerve(\we(\we(\mathcal M))^{[1]}) \rightleftarrows
\we(\mathcal M) \colon G \] which in turn induces a weak
equivalence of simplicial sets on the nerves, $U_\hoequiv \simeq
U_0$, which completes the proof.
\end{proof}

Now that we have proved that the mapping spaces and homotopy
categories agree for $V$ and for $L^H \mathcal M$, it remains to
show that they agree for $L^H \mathcal M$ and $L_C (L^H \mathcal
M)$.

\begin{theorem}
Let $\mathcal M$ be a model category.  For any $x$ and $y$ objects
of $L^H \mathcal M$, there is a weak equivalence of simplicial
sets
\[ \Map_{L^H \mathcal M}(x,y) \simeq \map_{L_C(L^H \mathcal M)}(x,y), \]
and there is an equivalence of categories $\pi_0 L^H \mathcal M
\approx Ho(L_C(L^H \mathcal M))$.
\end{theorem}

Note in particular that $x$ and $y$ are just objects of $\mathcal
M$, and that $\pi_0 L^H \mathcal M$ is equivalent to the homotopy
category $\Ho(\mathcal M)$.


\begin{proof}
Given the hammock localization $L^H \mathcal M$ of the model
category $\mathcal M$, we have the following composite map of
simplicial spaces:
\[ X= \nerve(L^H \mathcal M) \rightarrow X^f \rightarrow L_C(L^H \mathcal M).
\] Here, $X^f$ denotes a Reedy fibrant replacement of $X$.  This
composite is just Rezk's method for assigning the complete Segal
space $L_C(L^H \mathcal M)$ to the simplicial category $L^H
\mathcal M$.

On the left-hand side, the mapping spaces of $X= \nerve(L^H
\mathcal M)$ are precisely those of $L^H \mathcal M$, by the
definition of the nerve functor.  In the nerve, one of these
mapping spaces, say $\map_X(x,y)$ for some objects $x$ and $y$ of
$\mathcal M$, is given by the fiber over $(x,y)$ of the map $(d_1,
d_0) \colon X_1 \rightarrow X_0 \times X_0$.  Although these
mapping spaces can be defined for $X$, there is no reason that
they are homotopy invariant.  When we take a Reedy fibrant
replacement $X^f$ of $X$, however, this map becomes a fibration,
and hence this fiber is actually a homotopy fiber and so homotopy
invariant. For a general simplicial space, we cannot assume that
the mapping spaces of the Reedy fibrant replacement are equivalent
to the original ones. However, if the 0-space of the simplicial
space in question is discrete in degree zero, then the map above
is a fibration. Using an argument similar to the one in \cite[\S
5]{thesis}, we can find a Reedy fibrant replacement functor which
leaves the 0-space discrete. While the space in degree one might
be changed in this process of passing to $X^f_1$, it will still be
weakly equivalent $X_1$. In particular, the mapping spaces in
$X^f$ will be weakly equivalent to those in $X$.

Since the objects of $X^f$ are just the objects of $L^H \mathcal
M$, or the objects of $\mathcal M$, this equivalence of mapping
spaces gives us also an equivalence of homotopy categories.

The right-most map is the one defined by Rezk, $i_{X^f}: X^f
\rightarrow \widehat{X^f}$, which takes a Segal space to a
complete Segal space.  But, he defines this map in such a way that
it is in fact a Dwyer-Kan equivalence.  In other words, it induces
weak equivalences on mapping spaces and an equivalence of homotopy
categories. Thus, the composite map induces equivalences on
mapping spaces and an equivalence on homotopy categories.
\end{proof}

\section{A Characterization of Complete Segal Spaces Arising from
Simplicial Categories}

In this section, we give a thorough description of the weak
equivalence type of complete Segal spaces which occur as images of
Rezk's functor from the category of simplicial categories. We
consider several different cases, beginning with ones for which we
can use the classifying diagram construction, i.e., discrete
categories, and then proceed to more general simplicial
categories.

It should be noted that we are characterizing these complete Segal
spaces up to weak equivalence, and so the resulting descriptions
are of the homotopy type of the spaces in each simplicial degree.
For example, in the case of a discrete category, we describe the
corresponding complete Segal space in terms of the isomorphism
classes of objects, rather than in terms of individual objects, in
order to simplify the description.  One could just take all
objects, and generally get much larger spaces, if the more precise
description were needed for the complete Segal space corresponding
to a given category.

Furthermore, notice that determining the homotopy type of the
spaces in degrees zero and one are sufficient to determine the
homotopy type of all the spaces, since we are considering Segal
spaces.  Thus, we focus our attention on these spaces, adding in a
few comments about how to continue the process with the
higher-degree spaces.

\subsection{Case 1: $\mathcal C$ is a discrete groupoid}

If $\mathcal C=G$ is a group, then applying Rezk's classifying
diagram construction results in a complete Segal space equivalent
to $BG$, i.e., the constant simplicial space which is the
simplicial set $BG$ at each level. In particular, since all
morphisms are invertible, we obtain essentially no new information
at level 1 that we didn't have already at level 0.

\begin{example}
Let $G= \mathbb Z/2$.  Then $(NG)_0$ is just the nerve, or $B
\mathbb Z/2$.  Then $(NG)_1$ has two 0-simplices, given by the two
morphisms (elements) of $G$.  However, these two objects of
$G^{[1]}$ are isomorphic, and the automorphism group of either one
of them is $\mathbb Z/2$.  Thus, $(NG)_1$ is also equivalent to $B
\mathbb Z/2$.
\end{example}

If $\mathcal C$ has more than one object but only one isomorphism
class of objects, we get instead a simplicial space weakly
equivalent to the constant simplicial space which is $B\Aut(x)$ at
each level, for a representative object $x$.  If $\mathcal C$ has
more than one isomorphism class $\langle x\rangle$, then the
result will be weakly equivalent to the constant simplicial space
$\coprod_{\langle x\rangle}B\Aut(x)$.

\subsection{Case 2: $\mathcal C$ is a discrete category}

Since in the classifying diagram $N\mathcal C$, $(N\mathcal C)_0$
picks out the isomorphisms of $\mathcal C$ only, we still
essentially have $\coprod_{\langle x\rangle}B\Aut(x)$ at level 0.
However, if $\mathcal C$ is not a groupoid, then there is new
information at level 1. It instead looks like
\[ \coprod_{\langle x\rangle,\langle y\rangle} B\Aut(\coprod_
{\langle \alpha \rangle} \Hom (x,y)_\alpha) \] where the $\alpha$
index the isomorphism classes of elements of $\Hom (x,y)$.  The
subspace of $(N\mathcal C)_1$ corresponding to $\langle
x\rangle,\langle y\rangle$, denoted $(N \mathcal C)_1(x,y)$, fits
into a fibration
\[ \Hom (x,y) \rightarrow (N\mathcal C)_1(x,y) \rightarrow
B\Aut(x) \times B\Aut(y). \]


The space in dimension 2 is determined, then, by the spaces at
levels 0 and 1.  The subspace corresponding to isomorphism classes
of objects $\langle x \rangle, \langle y\rangle, \langle
z\rangle$, denoted $(N\mathcal C)_2(x,y,z)$, fits into a fibration
\[ \Hom (x,y) \times \Hom (y,z) \rightarrow (N\mathcal
C)_2(x,y,z) \rightarrow B\Aut(x) \times B\Aut(y) \times B\Aut(z).
\]  The whole space $(N \mathcal C)_2$, up to homotopy, looks like
\[ \coprod_{\langle x \rangle, \langle y \rangle, \langle z
\rangle} B\Aut \left( \coprod_{\langle \alpha \rangle, \langle
\beta \rangle} \Hom(x,y)_\alpha \times \Hom(y,z)_\beta \right). \]
We could describe each $(N \mathcal C)_n$ analogously.

\begin{example}
Let $\mathcal C$ denote the category with two objects and one
nontrivial morphism between them $(\cdot \rightarrow \cdot)$. If
$\{e\}$ denotes the trivial group, then $(N \mathcal C)_0 \simeq
B\{e\} \amalg B\{e\}$ and $(N\mathcal C)_1 \simeq B\{e\} \amalg
B\{e\} \amalg B\{e\}$. In particular, $N \mathcal C$ is not
equivalent to the classifying diagram of the trivial category with
one object and one morphism, which would be the constant
simplicial space $B\{e\}$.  However, note that the nerves of these
two categories are homotopy equivalent.  Thus, we can see that the
classifying diagram is more refined than the nerve in
distinguishing between these two categories.
\end{example}

\subsection{Case 3: $\mathcal C$ is a simplicial groupoid}

First, consider the case where we have a simplicial group $G$. Let
$G_n$ denote the group of $n$-simplices of $G$.  Then
\[  \hocolim_{\Deltaop}(\nerve(G_n)^t)=\nerve(G). \]

Let $L_C$ denote Rezk's completion functor which makes the nerve
into a complete Segal space.  We claim that
\[L_C(\hocolim_{\Deltaop}(\nerve(G_n)^t)) \simeq L_C(\hocolim_{\Deltaop}L_C (\nerve
(G_n)^t)).\]  We actually prove the more general statement that,
for any $X= \hocolim_{\Deltaop} X_n$,
\[L_C(\hocolim_{\Deltaop} X_n) \simeq L_C(\hocolim_{\Deltaop}L_C X_n).\]
To prove this claim, first note that we have Rezk's completion map
\[ i \colon \hocolim_{\Deltaop}X_n \rightarrow L_C(\hocolim_{\Deltaop}X_n) \]
which is a weak equivalence. Furthermore, since in $\css$ any
complete Segal space $Y$ is a local object and every object is
cofibrant, we have a weak equivalence of spaces
\[ \Map (L_C(\hocolim_{\Deltaop}X_n),Y) \simeq \Map(\hocolim_{\Deltaop}X_n,Y).
\]

So, for any complete Segal space $Y$, we have that
\[ \begin{aligned}
\Map(L_C \hocolim_{\Deltaop}(L_C X_n),Y) & \simeq
\Map(\hocolim_{\Deltaop}(L_C X_n),Y)
\\
& \simeq \holim_\Delta \Map (L_CX_n,Y) \\
& \simeq \holim_\Delta \Map (X_n,Y) \\
& \simeq \Map (\hocolim_{\Deltaop} X_n ,Y) \\
& \simeq \Map (L_C \hocolim_{\Deltaop}X_n,Y).
\end{aligned} \]

Note that the above calculation depends on the fact that,
\[ \Map(\hocolim_{\Deltaop} X_n,Y) \simeq \holim_{\bf \Delta}
\Map(X_n ,Y), \] which follows from working levelwise on
simplicial sets.

Then, since $G_n$ is a discrete group, completing its nerve is the
same as taking the classifying diagram $NG_n$ which, by case 1, is
weakly equivalent to the constant simplicial space $BG_n$, denoted
here $cBG_n$. Thus we have:
\[ \begin{aligned}
L_C(\nerve (G_n)) & \simeq L_C[\hocolim_{\Deltaop}(\nerve (G_n))] \\
& \simeq L_C[\hocolim_{\Deltaop}(L_C(\nerve (G_n)))] \\
& \simeq L_C[\hocolim_{\Deltaop} (cBG_n)] \\
& \simeq L_C(BG) \\
& \simeq BG.
\end{aligned} \]

So, we obtain a simplicial space weakly equivalent to the constant
simplicial space with $BG$ at each level. (Recall, however, that
$BG$ here is obtained by taking the diagonal of the simplicial
nerve, so it is not quite the identical case.)  If we have a
simplicial groupoid, rather than a simplicial group, we obtain the
analogous result, replacing $BG$ with
\[ \coprod_{\langle x\rangle}B\Aut(x). \]

\subsection{Case 4: $\mathcal C$ is a simplicial category with
every morphism invertible up to homotopy}

Alternatively stated, this case covers the situation in which
$\pi_0(\mathcal C)$ is a groupoid.

Recall that we have a model structure $\sco$ on the category of
categories with a fixed object set $\mathcal O$, in which the
cofibrant objects are retracts of free objects. So, taking a
cofibrant replacement of $\mathcal C$ in this model category
structure $\sco$ essentially gives a free replacement of $\mathcal
C$, denoted $F(\mathcal C)$, which is weakly equivalent to
$\mathcal C$.  (This cofibrant category can be obtained by taking
a simplicial resolution $F_* \mathcal C$ and then taking a
diagonal \cite[6.1]{dksimploc}.)  Now, taking the localization
with respect to all morphisms results in a simplicial groupoid.
So, we have Dwyer-Kan equivalences
\[ \xymatrix@1{F(\mathcal C)^{-1}F(\mathcal C) & F(\mathcal C)
\ar[l]_-{\simeq} \ar[r]^-\simeq & \mathcal C} \]

But, now $F(\mathcal C)^{-1}F(\mathcal C)$ is a simplicial
groupoid weakly equivalent to $\mathcal C$, so we have now reduced
this situation to case 3.

Note that, to write down a description of this complete Segal
space in terms of the original category $\mathcal C$, we need to
take isomorphism classes of objects in $\pi_0(\mathcal C)$, or
weak equivalence classes, as well as self-maps which are
invertible up to homotopy rather than strict automorphisms.  While
we will still use $\langle x \rangle$ to denote the equivalence
class of a given object, we will use $\Aut^h(x)$ to signify
homotopy automorphisms of $x$.  Thus, the complete Segal space
corresponding to $\mathcal C$ in this case essentially looks like
\[ \coprod_{\langle x \rangle} B\Aut^h(x) \] at each level.

\subsection{Case 5: $\mathcal C$ is any simplicial category}

First consider the subcategory of $\mathcal C$ containing all the
objects of $\mathcal C$ and only the morphisms of $\mathcal C$
which are invertible up to homotopy.  Apply case 5 to get a
complete Segal space, but take only the 0-space of it, to be the
0-space of the desired complete Segal space.

To find the 1-space, first recall the definition of the completion
functor as applied to a Segal space $W$:
\[ L_C(W) = L_{\css}(\diag([m] \mapsto (W^{E(m)})_n)). \]  Recall
further that $(W^{E(m)})_n= \Map(E(m) \times \Delta[n]^t, W)$.
Thus, the Segal space we obtain (before applying the functor
$L_{\css}$) looks like
\[ \Map(E(0) \times \Delta[0]^t,W) \Leftarrow \Map (E(1) \times
\Delta[1]^t,W) \Lleftarrow \Map(E(2) \times \Delta[2]^t,W) \cdots.
\]  If the Segal space $W$ is a fibrant replacement of
$\nerve(\mathcal C)$, then the space at level 1 consists of
diagrams
\[ \xymatrix{x \ar[r] \ar[d]_\simeq & y \ar[d]^\simeq \\
x' \ar[r] & y'} \] with the maps in the appropriate simplicial
level.

For simplicity, we restrict to a given pair of objects $x$ and
$y$, representing given equivalence classes.  Consider the
homotopy automorphisms of $x$ and $y$.  If they are not all
invertible, we take a cofibrant replacement and group completion
as in case 4. So, without loss of generality, assume that
$\Aut(x)$ and $\Aut(y)$ are simplicial groups.  Note that we have
\[ \Aut(x)=\hocolim_{\Deltaop} \Aut(x)_n \]
and
\[ \Aut(y)=\hocolim_{\Deltaop} \Aut(y)_n. \]

Now look at
\[ \Map (x,y)= \hocolim_{\Deltaop} \Map (x,y)_n. \]
Consider for each $n \geq 0$ the discrete category $\mathcal C
(x,y)_n$ which has as objects $\Map (x,y)_n$ and as morphisms
pairs $(\alpha ,\beta)$ of automorphisms in $\Aut(x)_n \times
Aut(y)_n$ making a commutative square
\[ \xymatrix{x \ar[r]^f \ar[d]_\alpha & y \ar[d]^\beta
\\
x \ar[r]^{f'} & y} \] with $f,f' \in \Map (x,y)_n$.

Thus, the 1-space that we are interested in is also the 1-space of
the complete Segal space given by
\[ L_{\css}(\hocolim_{\Deltaop} (\nerve(\mathcal C(x,y)_n))). \]
Using a straightforward argument about localization functors
similar to the one in case 3 (which can be found in
\cite[4.1]{simpmon}), we can also apply the functor $L_{\css}$ on
the inside to get an equivalent simplicial space
\[ L_{\css}(\hocolim_{\Deltaop} (L_{\css}\nerve(\mathcal C(x,y)_n))). \]
But, since $\mathcal C(x,y)_n$ is a discrete category, this space
is just
\[ L_{\css} \hocolim_{\Deltaop} (N\mathcal C(x,y)_n) \simeq \hocolim_{\Deltaop} (N \mathcal C(x,y)_n). \]
Now, we restrict to the 1-space here, which is
\[ \hocolim_{\Deltaop}\left( \coprod_{\langle x \rangle, \langle
y \rangle} B\Aut \left(\coprod_{\langle \alpha \rangle} \Map_n
(x,y)_\alpha \right) \right) \simeq \coprod_{\langle x \rangle,
\langle y \rangle} B\Aut \left(\coprod_{\langle \alpha \rangle}
\Map(x,y)_\alpha \right). \]

As with the previous case, we can then go back and weaken to
homotopy automorphisms and equivalence classes of objects to
consider categories before taking a group completion, so our space
looks like
\[ \coprod_{\langle x \rangle, \langle y
\rangle} B\Aut^h \left(\coprod_{\langle \alpha \rangle}
\Map(x,y)_\alpha \right). \]

We could then obtain the 2-space of our complete Segal space by
considering categories $\mathcal C(x,y,z)_n$ defined similarly,
and the description of the 2-space of the classifying diagram of a
discrete category as given in case 2.

We can summarize these results in the following theorem.  For an
object $x$ of a simplicial category $\mathcal C$, let $\langle
x\rangle$ denote the weak equivalence class of $x$ in $\mathcal
C$, and for a morphism $\alpha \colon x \rightarrow y$, let
$\langle\alpha\rangle$ denote the weak equivalence class of
$\alpha$ in the morphism category $\mathcal C^{[1]}$. Let
$\Aut^h(x)$ denote the space of self-maps of $x$ which are
invertible in $\pi_0\mathcal C$.

\begin{theorem}
Let $\mathcal C$ be a simplicial category.  The complete Segal
space corresponding to $\mathcal C$ has the form
\[ \coprod_{\langle x\rangle} B\Aut^h(x) \Leftarrow \coprod_{\langle x\rangle,\langle y\rangle} B\Aut^h
\left(\coprod_{\langle \alpha \rangle} \Map(x,y)_\alpha \right)
\Lleftarrow \cdots.
\]
\end{theorem}

\end{document}